\documentclass{amsart}
\usepackage{amsmath,amssymb,amsthm,amscd}
\usepackage{graphicx}
\usepackage{hyperref}

\title{Sparse stable packings of hard discs in a box}
\author{Matthew Kahle}
\address{Institute for Advanced Study}
\email{mkahle@math.ias.edu}
\thanks{Supported in part by Stanford's NSF-RTG grant in geometry and topology}

\date{\today}

\begin{document}

\newtheorem{theorem}{Theorem}[section]
\newtheorem{lemma}[theorem]{Lemma}
\newtheorem{corollary}[theorem]{Corollary}
\newtheorem{proposition}[theorem]{Proposition}
\theoremstyle{definition}
\newtheorem{definition}[theorem]{Definition}
\newtheorem{example}[theorem]{Example}
\newtheorem{claim}[theorem]{Claim}

\theoremstyle{remark}
\newtheorem{remark}[theorem]{Remark}
\newcommand{\be}{\begin{equation}}
\newcommand{\ee}{\end{equation}}
\newcommand{\prob}{\mbox{\bf P}}
\newcommand{\Z}{\mathbb{Z}}
\newcommand{\R}{\mathbb{R}}
\newcommand{\Q}{\mathbb{Q}}
\newcommand{\config}{\mbox{Config}}



\begin{abstract} We construct stable, or locally jammed, configurations of $n$ overlapping discs of radius $r$ in a unit square, with $r = O(1/n)$.  This is best possible, up to a constant factor. A consequence is that the Metropolis algorithm, a well-studied Markov chain on the hard discs model, is not ergodic in this range of parameters as is commonly assumed.
\end{abstract}

\maketitle


\section{Introduction}

Sphere packing problems have been studied for at least the past few centuries, since Kepler and Newton. The densest packing of spheres in three dimensions is easy to guess and must have been known for quite a while, and was conjectured by Kepler, but the fact was only proved recently by Thomas Hales \cite{Hales}. (Under the assumption that the centers of the spheres lie on a lattice, the problem is much easier, and the answer was already known to Gauss.) Hales' proof famously relies on a difficult computer calculation.

Even the two-dimensional case of sphere packing, the densest packing of circles in the plane, is nontrivial. This problem was solved by L. Fejes T\'oth only in 1940. Since then, there has also been a wide literature on densest packing of spheres in a bounded convex region, such as a circle, triangle, or square \cite{Graham1,Graham2,Graham3, Graham4}. The case of a region with a boundary seems much harder in general.

We use the words ``disc packing'' in this article to distinguish from the area of circle packings and discrete conformal geometry \cite{circle_packing}. As an interesting historical note on disc packings, in 1964 B\"or\"oczky used a subtle construction to disprove a conjecture of Fejes T\'oth, that a stable arrangement of discs in the plane, one where each disc is held in place by its neighbors, must have positive density \cite{Boro}. B\"or\"oczky's counterexample is at the heart of what we do here.

\begin{figure}
\includegraphics[width=3.5in]{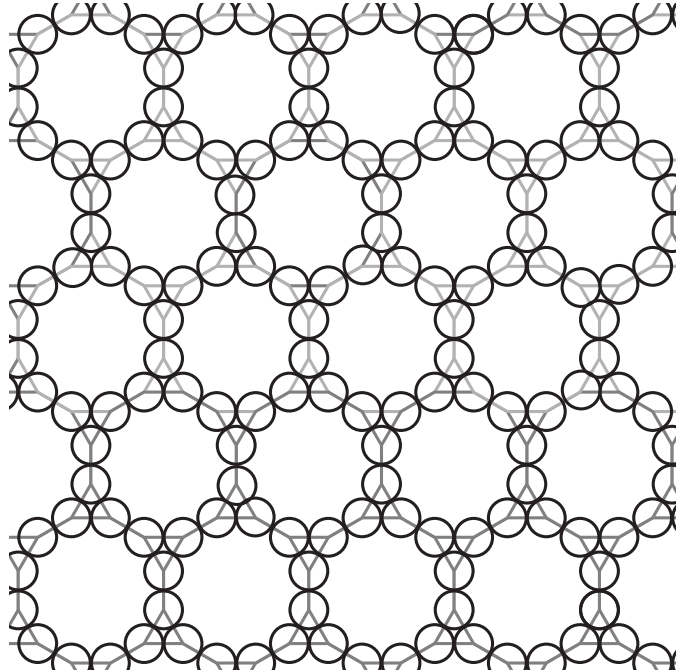}
\caption{A tiling by regular dodecagons and equilateral triangles, and placing suitably scaled discs at all of the vertices, results in a stable configuration in the plane.}
\label{fig:sparse-stable}
\end{figure}

\begin{figure}
\includegraphics[width=3.5in]{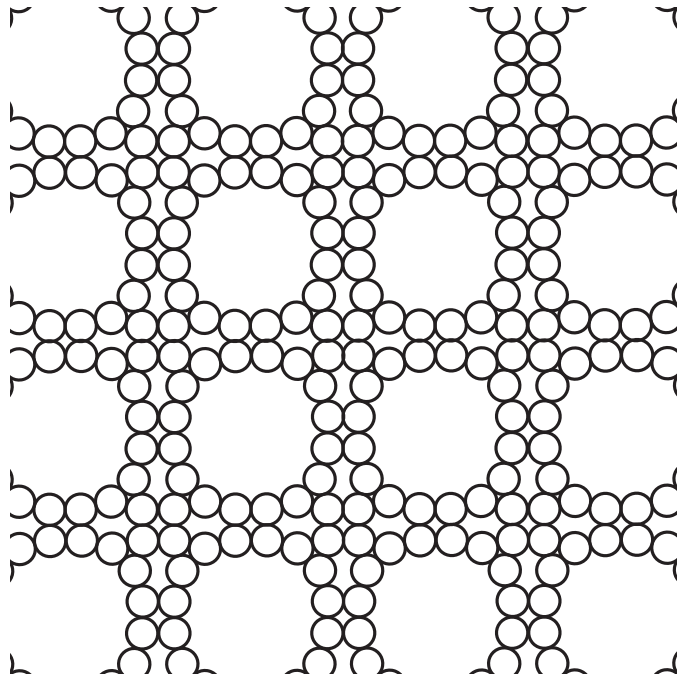}
\caption{Another planar stable configuration. This also gives examples of stable configurations in a square, for a sequence of increasing $n$ and $ r  = \Omega( 1/ \sqrt{n})$.}
\label{fig:rigid5-planar}
\end{figure}

Regular $12$-gons and equilateral triangles are well known to tile the plane vertex-transitively \cite{Tilings}. If suitably scaled discs are placed at the vertices of the polygons, one obtains a stable configuration, as in Figure \ref{fig:sparse-stable}. Apparently some have conjectured that this configuration has minimal density \cite{Penguin}. One computes that these discs asymptotically cover a fraction of $(7 \sqrt{3} -12) \pi \approx 0.390675$ of the plane. We note another planar stable configuration in Figure \ref{fig:rigid5-planar}, which we discovered while preparing this article. This planar configuration has a slightly greater density of $0.4569$, but larger empty regions. We note that this configuration easily embeds isometrically in a square.

Besides the discrete geometry background, we are also motivated by a problem in statistical physics, and in fact this was our original motivation. Hard discs in a box is one of the most basic mathematical models for matter. Here one considers the space of all possible configurations of $n$ non-overlapping discs of radius $r$ in a unit square. This space has a natural probability measure, which we will not discuss here. But conceptually one could think of noninteracting particles, or particles where attractive forces are so much weaker than repulsive forces that the repulsive forces dominates, and so in the limit making a constraint to non overlapping particles makes sense. See the background and references in Section 4 of \cite {Persi_Markov}. It is of theoretical interest to try to understand basic properties of phase diagrams, even in this simple mathematical model.

Countless hours of computer simulations have been run with hundreds of discs, and various densities of $r$. (For example, see Chapter 3 of \cite{Krauth}). It has been observed that there is a phase transition when $\pi r^2 n \approx 0.70$. In particular if $\pi r^2 n  \ll 0.70$ then the discs seem to be behave more ``liquid,'' and if $\pi r^2 n \gg 0.70$, they are more ``solid,'' according to various statistics measured, indicating action at a distance, and so on. To sample a point in $\config(n;r)$ uniformly, the Metropolis algorithm and its variants are used.

The most basic implementation of Metropolis is to start with a particular configuration of discs, then iterates the following algorithm: choose one of the $n$ discs uniformly randomly, and then perturb it randomly, by moving it uniformly within a smaller disc of radius $\epsilon$. If the proposed new position of the disc does not go out of the box, and does not overlap with any other discs, then accept the move. Otherwise, reject it. Repeat this process as necessary.

In order to make the ``as necessary'' precise, Diaconis, Lebeau, and Michel recently established useful bounds on the mixing time of the Metropolis algorithm in general, and then applied their results to the specific case of hard discs in a box \cite{Persi_geometric}. In order to meaningfully discuss mixing time questions, they must require that the discs be sufficiently small that one can move from every configuration to every other configuration, via small moves by one disc at a time (i.e. that the Metropolis Markov chain is irreducible), and they give an upper bound on the disc radius that guarantees that this is the case. Let $\config(n;r)$ denote the space of all possible configurations of $n$ non-overlapping discs of radius $r$ in a unit square. A stable configuration represents an isolated point in $\config(n;r)$.

\begin{theorem}[Diaconis, Lebeau, and Michel] \label{irreducible} There exists a constant $\alpha$ such that the Metropolis Markov chain on $\config(n;r)$ is irreducible whenever $r \le \alpha / n$. In particular, there are no stable configurations in this range.
\end{theorem}

One might wonder if this result can be extended to larger $r$. Is it still true, for example, if $r =O( \sqrt{n})$? However, just as in Fejer T\'oth's boundary-free version of this problem, the answer is no. The result of Diaconis, Lebeau, and Michel is is essentially best possible.

\begin{theorem} \label{theorem-stable} There exists a constant $\beta$ such that for arbitrarily large $n$ and $r \ge \beta / n$, there exist stable configurations of $n$ discs of radius $r$.
\end{theorem}

This is our main result. We prove Theorem \ref{theorem-stable} in the next section.

\section{Sparse stable configurations} \label{section:stable} A positioning of $n$ nonoverlapping discs of radius $r$ in the unit square $[0,1]^2$ is a {\it configuration}. We say that a configuration is {\it stable} if each disc is held in place by its neighbors and the walls. Our main result is the following.

\begin{figure}
\includegraphics[width=3.5in]{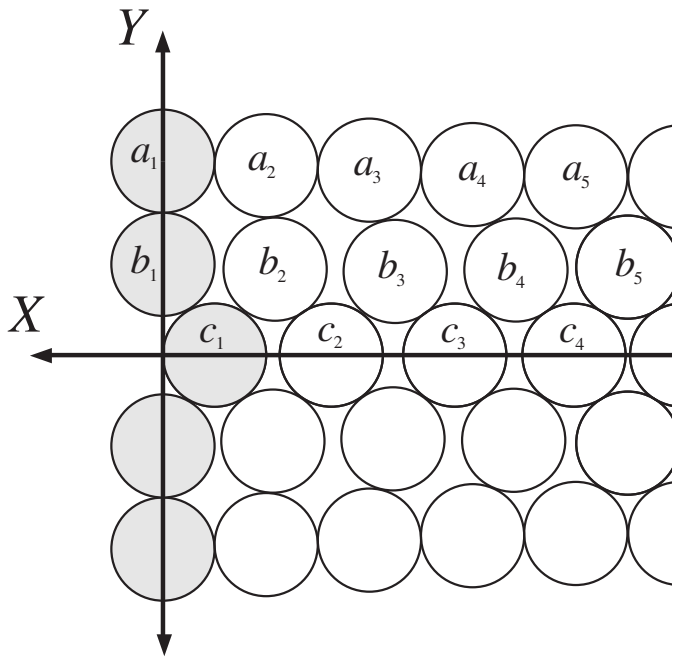}
\caption{B\"or\"oczky's bridge. (After Figure 9.2, p. 282, in \cite{Pach}.)}
\label{fig:bridge1}
\end{figure}

We first describe a construction of a one-way infinite ``bridge'' in the plane, due to B\"or\"oczky \cite{Boro}. We follow closely Pach and Sharir's description in \cite{Pach}. For convenience of notation, all of the discs for now are unit radius; they can be rescaled as necessary later to fit in a box. First, the construction is symmetric about the $x$-axis, so we only describe discs with their centers on our above the axis. Let $a_1 = (0, 2 + \sqrt{3}) $, $b_1 = (0, \sqrt{3})$, $c_1=(1,0)$, as in Figure \ref{fig:bridge1}. Let $f(x)$  be strictly convex function defined for $x \ge 0$, with $f(0) = 2 + \sqrt{3}$ and $\lim_{x \to \infty} f(x) = 2 \sqrt{3}.$ Let $C$ be the curve $C = \{ (x,f(x)) \mid x\ge 0 \}$, and let $a_1, a_2, a_3, \ldots$, the unique sequence of distinct points on $C$ such that $d(a_i, a_{i+1}) = 2 $ for $i = 1, 2, \ldots$.

Now set $b_2$ to be the unique point to the right of $a_2$ such that $d(b_2, a_2)=2$, and $d(b_2, c_1)=1$. Then set $c_2$ to be the point on the $x$-axis to the right of $b_2$ with $d(c_2,b_2)=1$. Inductively define $b_{i+1}$ to be the point to the right of $a_i$ such that $d(b_{i+1},a_{i+1}) = 2$ and $d(b_{i+1},c_i)=2$. Similarly, set $c_{i+1}$ to be the point on the $x$-axis to the right of $b_{i+1}$ such that $d(c_{i+1},b_{i+1})=2$.

\begin{figure}
\includegraphics[width=4in]{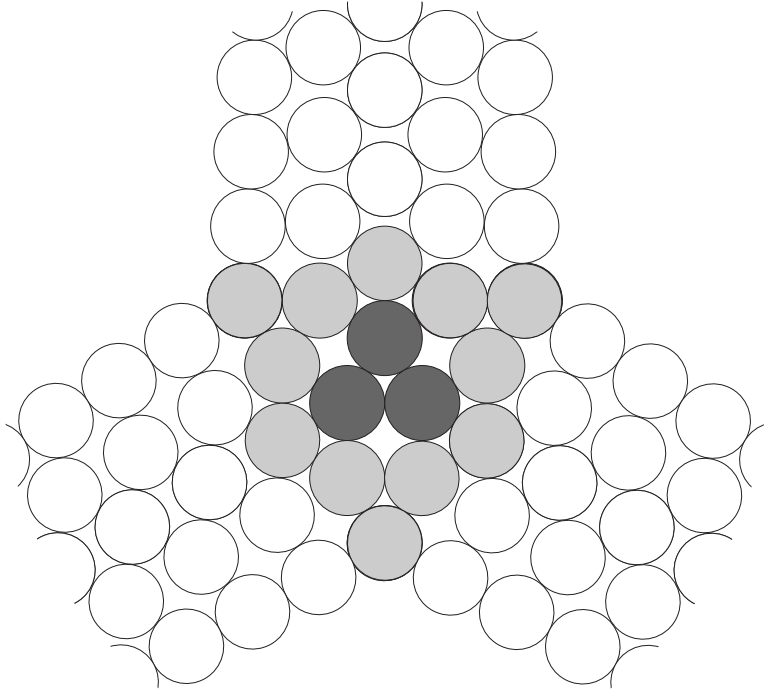}
\caption{B\"or\"oczky's junction. The dark shaded discs are stable, and three copies of the bridge will overlap along the lightly shaded discs. This gives a zero-density stable configuration in the plane.}
\label{fig:Bjunction}
\end{figure}

\begin{figure}
\includegraphics[width=5.5in]{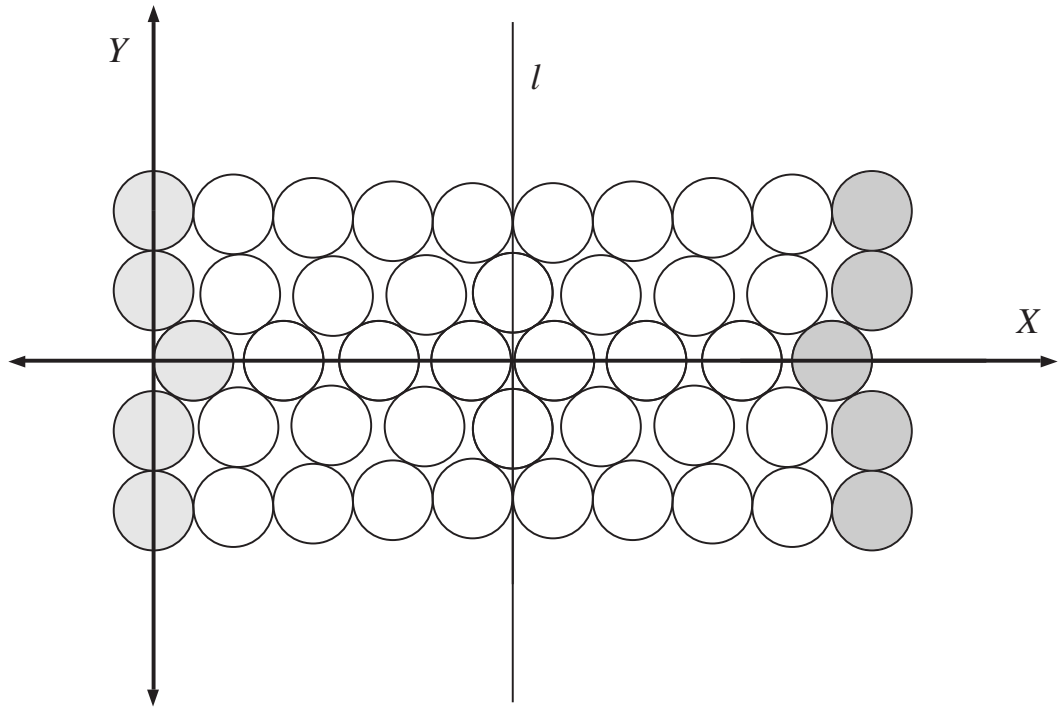}
\caption{Pach and Sharir's symmetric version of B\"or\"oczky's bridge. \vspace{.5in}}
\label{fig:bridge2}
\end{figure}

The beginning of this construction is illustrated in the left side of Figure \ref{fig:bridge2}. The details that the construction results in a well-defined, stable configuration (except along the four discs on the left) can be found in B\"or\"oczky's original paper \cite{Boro}. He needed a second idea to get a stable configuration and disprove Fejer Toth's conjecture. He constructed a junction where three of these bridges could meet and hold all the loose discs on the end in place. B\"or\"oczky's junction is shown in Figure \ref{fig:Bjunction}.

We first use a slight modification of B\"or\"oczky's bridge construction, due to Pach and Sharir \cite{Pach}. For $\epsilon > 0$, replace $f(x)$ above by the strictly convex function $$f_{\epsilon}(x):=(1+\epsilon)f(x) - \epsilon f(0).$$ We have $f_{\epsilon} (0) = f(0) = 2 + \sqrt{3}$, and $\lim_{x \to \infty} f_{\epsilon}(x) < 2\sqrt{3}$, so the sequence no longer continues indefinitely. By the intermediate value theorem, by varying $\epsilon$ we can insure that for any arbitrarily large $N$, some $b_N$ is the last well-defined point, and its $x$-coordinate is one more than the $x$-coordinate of $a_{N}$.

Then let $l$ be the vertical line through $b_N$, as in Figure \ref{fig:bridg2e}, and complete the configuration using $l$ as a line of symmetry, as illustrated. This gives an arbitrarily long symmetric bridge.

\begin{figure}
\includegraphics[width=2.75in]{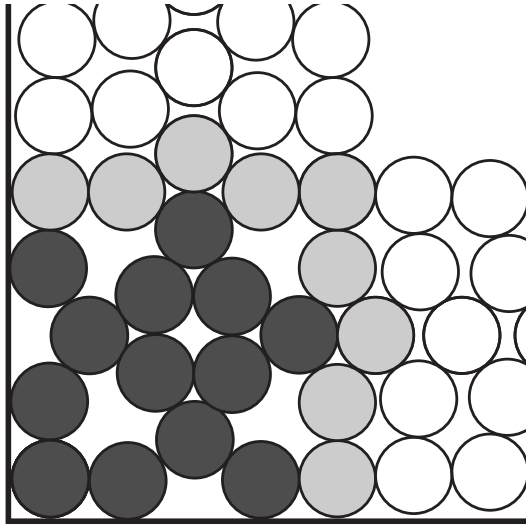}
\caption{The corner junction.}
\label{fig:junction}
\end{figure}

The second piece of our construction for the square is the corner ``junction'' piece in Figure \ref{fig:junction}. The construction is easily verified to exist. If the center of the bottom left disc has coordinates $(0,0)$, then the remaining five discs in the bottom left quarter of the square have coordinates: $$(0,2), (2,0), (2+\sqrt{3}, 1), (1, 2+\sqrt{3}), \mbox{ and }(1+\sqrt{3},1+\sqrt{3}).$$

\begin{figure}
\includegraphics[width=5in]{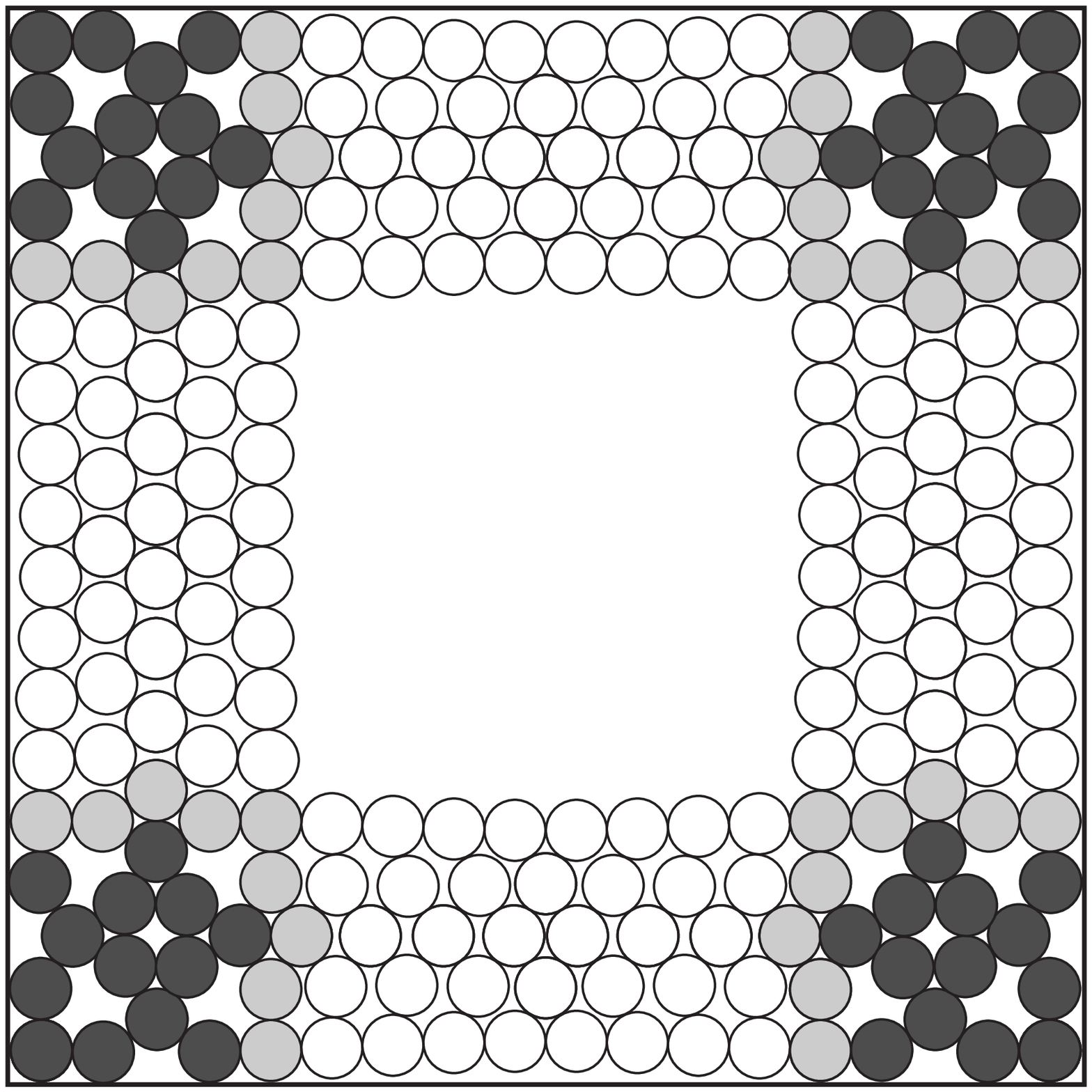}
\caption{The main construction. The white bridges can be stretched arbitrarily long, so there are stable configurations in the square with arbitrarily large $n$ and $r=O(1/n)$.}
\label{fig:rigid}
\end{figure}

Then by taking four copies each of the junction and bridge, and arranging them so that they overlap along their shaded outermost disks, in Figure \ref{fig:rigid}, the resulting configuration is stable.

\section{Future directions} \label{future}

The configurations constructed in this article are stable, meaning that no single disc can move. This is important to know, particularly from the point of view of probability and statistical physics. But is it possible for several discs to move simultaneously? Given that one can continuously deform the strictly convex function $f(x)$ in the main construction in Section \ref{section:stable}, and if one is careful preserve the combinatorial structure (tangency graph), it would seem that there is a small amount of flexibility and that one can slightly perturb many disks at once. How much freedom is there, though? In particular, is it possible to freely permute the discs?

The set of all configurations of $n$ hard discs in a unit box, $\config(n;r)$, inherits a subspace topology from $\R^{2n}$. Persi Diaconis points out in his survey article \cite{Persi_Markov} that, ``very, very little is known about the topology of these spaces.'' The question above about moving multiple discs simultaneously is most naturally phrased in terms of path components of $\config(n;r)$. Considering $\config(n;r)$ as a topological space allows us to ask many new questions, as well. For example, we could consider path components under various restrictions (Lipshitz, smooth, etc.).  If $r$ is sufficiently small, it would seem that $\config(n;r)$ is not only path connected, but homotopy equivalent to the classical configuration space $\config(n, \R^2)$, a space well-studied in algebraic geometry and algebraic topology \cite{Fulton}. What is the largest $R = R(n)$ such that $\config(n;r)$ is homotopy equivalent to $\config(n, \R^2)$ for $r \le R$?

\begin{figure}
\includegraphics[width=2.5in]{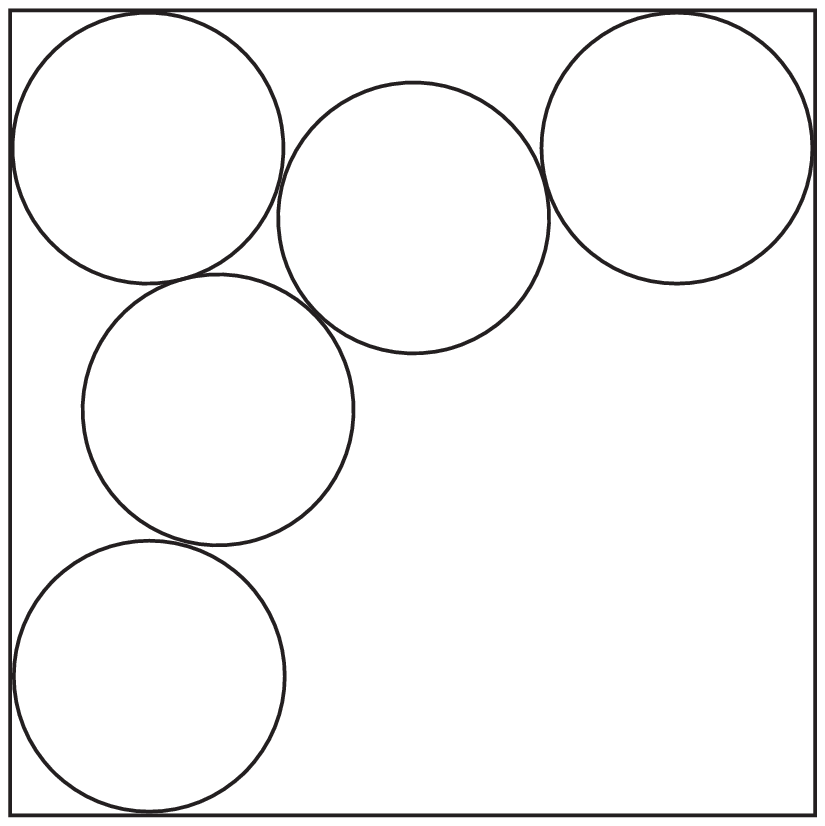}
\caption{A stable configuration of $5$ discs that illustrated non-monotonicity of the property ``path-connectedness.'' This configuration is an isolated point in $\config(n;r)$, but for slightly larger radius, computer experiments suggest that $\config(5;r)$ is path connected.}
\label{fig:rigid5}
\end{figure}

One might guess that once $r$ is sufficiently small that $\config(n;r)$ is path connected, that if $r' < r$ then $\config(n;r')$ is path connected as well. However, this is not the case. For example, the configuration in Figure \ref{fig:rigid5} is stable, and in fact, it can be shown, that it is a single point component in $\config(n;r)$. But for slightly larger radius, it would seem that $\config(n;r)$ is path connected.

On the subject of non-monotonicity, we note the work of Brightwell, H{\"a}ggstr{\"o}m, and Winkler , where hardcore and Widom-Rowlinson models are studied on graphs. In \cite{Winkler-nonmonotone} the property of having a unique Gibbs measure, perhaps somewhat analogous to $\config(n;r)$ being path connected, is shown to be non-monotone in the underlying parameter.

The topological questions may be interesting for their own sake, but they may also have implications for understanding the statistical physics. For example, although we have showed that Metropolis is not irreducible on this problem, even for reasonably small discs, some might complain that the stable configurations only comprise a set of measure zero in $\config(n;r)$. Suppose that one could show that $\config(n;r)$ is disconnected. Then one can use a compactness argument to show that $\config(n;r')$ is still disconnected for some $r' < r$. In fact for some $r' < r$, $$ |\pi_0(\config(n;r)) | <  |\pi_0(\config(n;r')) | ,$$ where $\pi_0(X)$ denotes the set of path components of $X$. There is a natural inclusion map  $i: \config(n;r) \hookrightarrow \config(n;r'),$ and every element of $i(\pi_o(\config(n,r)) $ has positive measure in $\config(n;r')$.

For example, in Figure \ref{fig:rigid5}, one can slightly shrink the discs, and there are still at least $4 \times 5! +1 = 481$ path components, each with small positive measure.

Other discrete geometry questions naturally arise. Can an analogue of the construction in this article be carried out in higher dimensional cubes? What about for an arbitrary smoothly bounded convex body in the plane? (The result of Diaconis, et. al \cite{Persi_geometric} is much more general than we stated earlier. In particular, Theorem \ref{irreducible} holds for general convex domains with reasonably nice boundaries.)

The most important question in this area, from the point of view of statistical physics, is to give a mathematical explanation for the phase transition observed experimentally to occur when $\pi r^2 n \approx 0.70.$ To prove that any phase transition occurs at all in the hard spheres model seems to be an open problem, although we note \cite{Bowen}, where a a phase transition is proved to exist, in a similar model where the particles are nonconvex.

\section*{Acknowledgements} The author thanks Jackson Gorham for pointing out the stable configuration in Figure \ref{fig:rigid5}. Jackson is an undergraduate researcher at Stanford, who discovered the configuration with Matlab code, by flowing to locally ``minimal energy'' configurations. Especially, thanks to Persi Diaconis for pointing out this problem and for suggestions which greatly improved the article.


\begin{thebibliography}{10}

\bibitem{Graham2}
David~W. Boll, Jerry Donovan, Ronald~L. Graham, and Boris~D. Lubachevsky.
\newblock Improving dense packings of equal disks in a square.
\newblock {\em Electron. J. Combin.}, 7:Research Paper 46, 9 pp. (electronic),
  2000.

\bibitem{Boro}
K.~B{\"o}r{\"o}czky.
\newblock \"{U}ber stabile {K}reis- und {K}ugelsysteme.
\newblock {\em Ann. Univ. Sci. Budapest. E\"otv\"os Sect. Math.}, 7:79--82,
  1964.

\bibitem{Bowen}
L.~Bowen, R.~Lyons, C.~Radin, and P.~Winkler.
\newblock A solidification phenomenon in random packings.
\newblock {\em SIAM J. Math. Anal.}, 38(4):1075--1089 (electronic), 2006.

\bibitem{Winkler-nonmonotone}
Graham~R. Brightwell, Olle H{\"a}ggstr{\"o}m, and Peter Winkler.
\newblock Nonmonotonic behavior in hard-core and {W}idom-{R}owlinson models.
\newblock {\em J. Statist. Phys.}, 94(3-4):415--435, 1999.

\bibitem{Persi_Markov}
Persi Diaconis.
\newblock The {M}arkov chain {M}onte {C}arlo revolution.
\newblock {\em Bull. Amer. Math. Soc. (N.S.)}, 46(2):179--205, 2009.

\bibitem{Persi_geometric}
Persi Diaconis, Gilles Lebeau, and Laurent Michel.
\newblock Geometric analysis for the {M}etropolis algorithm on {L}ipshitz
  domains.
\newblock {\em Technical report, Department of Statistics, Stanford
  University}, 2008.
\newblock preprint.

\bibitem{Fulton}
William Fulton and Robert MacPherson.
\newblock A compactification of configuration spaces.
\newblock {\em Ann. of Math. (2)}, 139(1):183--225, 1994.

\bibitem{Graham1}
R.~L. Graham and B.~D. Lubachevsky.
\newblock Repeated patterns of dense packings of equal disks in a square.
\newblock {\em Electron. J. Combin.}, 3(1):Research Paper 16, approx.\ 17 pp.\
  (electronic), 1996.

\bibitem{Tilings}
Branko Gr{\"u}nbaum and G.~C. Shephard.
\newblock {\em Tilings and patterns}.
\newblock W. H. Freeman and Company, New York, 1987.

\bibitem{Hales}
Thomas~C. Hales.
\newblock A proof of the {K}epler conjecture.
\newblock {\em Ann. of Math. (2)}, 162(3):1065--1185, 2005.

\bibitem{Krauth}
Werner Krauth.
\newblock {\em Statistical mechanics}.
\newblock Oxford Master Series in Physics. Oxford University Press, Oxford,
  2006.
\newblock Algorithms and computations, Oxford Master Series in Statistical
  Computational, and Theoretical Physics.

\bibitem{Graham3}
B.~D. Lubachevsky and R.~L. Graham.
\newblock Curved hexagonal packings of equal disks in a circle.
\newblock {\em Discrete Comput. Geom.}, 18(2):179--194, 1997.

\bibitem{Graham4}
Boris~D. Lubachevsky, Ron~L. Graham, and Frank~H. Stillinger.
\newblock Patterns and structures in disk packings.
\newblock {\em Period. Math. Hungar.}, 34(1-2):123--142, 1997.
\newblock 3rd Geometry Festival: an International Conference on Packings,
  Coverings and Tilings (Budapest, 1996).

\bibitem{Pach}
J{\'a}nos Pach and Micha Sharir.
\newblock {\em Combinatorial geometry and its algorithmic applications}, volume
  152 of {\em Mathematical Surveys and Monographs}.
\newblock American Mathematical Society, Providence, RI, 2009.
\newblock The Alcal{\'a} lectures.

\bibitem{circle_packing}
Kenneth Stephenson.
\newblock {\em Introduction to circle packing}.
\newblock Cambridge University Press, Cambridge, 2005.
\newblock The theory of discrete analytic functions.

\bibitem{Penguin}
David Wells.
\newblock {\em The {P}enguin dictionary of curious and interesting geometry}.
\newblock Penguin Books, New York, 1991.

\end{thebibliography}
\end{document}